\documentclass[final]{dmtcs-episciences}

\usepackage[utf8]{inputenc}
\usepackage{subfigure}

\usepackage[round]{natbib}

\usepackage[english]{babel}
\usepackage{amssymb}
\usepackage{amsmath}
\usepackage{latexsym}
\usepackage{amsthm}
\usepackage[pdftex]{graphicx}
\usepackage{tikz} 
\usepackage{hyperref}
\usepackage{pgfplots}  
\usepackage{mathtools}
\pgfplotsset{width=6.6cm,compat=1.7}

\usetikzlibrary{cd}

\usetikzlibrary{decorations.pathreplacing}

\DeclareMathOperator{\fp}{fp}

\DeclareMathOperator*{\cyc}{cyc}

\newcommand{\seqnum}[1]{\href{http://oeis.org/#1}{\underline{#1}}}

 \theoremstyle{plain}
 \newtheorem{thm}{Theorem}[section]

 \theoremstyle{definition}

 \theoremstyle{remark}
 
\author[Marilena Barnabei et al.]{Marilena Barnabei\affiliationmark{1}
  \and Niccol\`o Castronuovo\affiliationmark{2}
  \and Matteo Silimbani\affiliationmark{3}}
\title[Hertzsprung patterns on involutions]{Hertzsprung patterns on involutions}

\affiliation{
  Universit\`a di Bologna, Bologna, Italy \\
  Liceo ``A. Einstein'', Rimini,  Italy \\
  Istituto Comprensivo ``E. Rosetti'', Forlimpopoli, Italy }

\begin{document}
\publicationdata{vol. 27:1, Permutation Patterns 2024}{2025}{15}{10.46298/dmtcs.14897}{2024-12-05; 2024-12-05; 2025-05-07}{2025-09-03}
\maketitle

\begin{abstract}
\vspace{0.25cm}
Hertzsprung patterns, recently introduced by Anders Claesson, are subsequences of a permutation contiguous in both positions and values, and can be seen as a subclass of bivincular patterns.

This paper investigates Hertzsprung patterns within involutions, where additional structural constraints introduce new challenges. We present a general formula for enumerating occurrences of these patterns in involutions.

We also analyze specific cases to derive the distribution of all Hertzsprung patterns of lengths two and three.
\end{abstract}

In memory of Flavio Bonetti (1948 - 2023).

\section{Introduction}
The study of consecutive patterns in permutations, initiated in \cite{Babson2000},   has garnered significant attention in combinatorics, motivated by its connections to diverse fields such as algebra, geometry, and computer science.

Consecutive patterns capture the relative order of elements in contiguous subsequences of permutations, providing a refined lens through which to analyze structural properties. 

The concept of Hertzsprung pattern - an addition to the landscape of permutation patterns recently introduced by \cite{Claesson_Hertz} - focuses on identifying subsequences contiguous not only in positions but also in values. Hertzsprung patterns can be seen as a special instance of  bivincular patterns, first introduced in \cite{BOUSQUETMELOU2010884}. Bivincular patterns generalize classical and vincular patterns by incorporating dependencies between non-adjacent positions.
In \cite{Claesson_Hertz}, the author determines the joint distribution of occurrences of any set of (incomparable) Hertzprung patterns using a variation of the Goulden-Jackson cluster method (\cite{Goulden}).

Drawing inspiration from the work of Claesson, in this paper we explore the enumeration of Hertzsprung patterns within the set of involutions - permutations that are their own inverses. This study presents some additional complexities than in general permutations, due to the structural constraints imposed by the self-inverse nature of involutions. 

Using the combinatorial theory of generating functions in form of continued fractions and a variant of the Goulden-Jackson method, in Section \ref{general_results} we provide a general formula for the enumeration of occurrences of Hertzsprung patterns in involutions. 
In the subsequent sections we analyze various specific cases that allow us to determine the distribution of all Hertzsprung patterns of length two and three.
The last section is devoted to some considerations concerning\\ Wilf-equivalences among Hertzprung patterns of length two or three.

\section{Preliminaries} 

We will use the following terminology, as in \cite{Claesson_Hertz}. Let $\pi,$ $\beta$ be words. The word $\beta$ is said to be a \emph{factor} of $\pi$ if there are words $\alpha$ and $\gamma$ such that $\pi=\alpha \beta \gamma$. 
We will denote by $|\pi|$ the length of the word $\pi$.
Let $\mathcal S_k$ be the symmetric group over $\{1,\ldots,k\}$. Let $\tau \in \mathcal S_k,$ with $k\geq 2,$ and let $\pi$ be a word over the set of positive integers without repetitions.  We call $\tau$ a \emph{Hertzsprung pattern} (\emph{H-pattern}, for short) of $\pi$ if there are an integer $c$ and a factor $\beta=b_1b_2\ldots b_k$ of $\pi$ such that $b_i-\tau_i=c$ for each $i\in \{1,2,\ldots,k\}.$ In this case we say that $\beta$ is an \emph{occurrence} of the H-pattern $\tau$ and write $\tau \leq \pi.$ The number of occurrences of the H-pattern $\tau$ in $\pi$ will be denoted by $\tau(\pi).$ For
an example, 231 is a H-pattern of $\pi = 897235641$; indeed, $897$ and $564$ are the only occurrences of the H-pattern 231 in $\pi,$  hence $231(\pi)=2.$ 
Notice that a Hertzprung pattern $\tau_1\tau_2\ldots \tau_k$ can be seen as the bivincular pattern $\overline{\underline{\tau_1\tau_2\ldots \tau_k}}$ (as defined in \cite[pg. 13]{Ki}).

A set of H-patterns will be called an \emph{antichain} when none of its elements is a H-pattern of the others. 
From now on, $T$ will denote an antichain of H-patterns.
A $T$-\emph{marked word} is a pair $(\pi,M),$ where $\pi$ is a word over the set of positive integers without repeated symbols and $M$ is a (possibly empty) subset of all occurrences in $\pi$ of patterns from $T$. Every element of $M$ is called a $T$-\emph{marked occurrence}. 

As an example consider $T=\{2134,123\}$ and let $\pi=1\,3\,2\,4\,5\,6\,8\,7\,9\,10.$
Observe that $\pi$ contains two occurrences of the H-pattern 2134, namely, $3\,2\,4\,5,$ and $8\,7\,9\,10$, and a single occurrence of the H-pattern 123, namely, $4\,5\,6.$ In this case, the pairs $(\pi,\{4\,5\,6,3\,2\,4\,5,8\,7\,9\,10\}),$ $(\pi,\{4\,5\,6,8\,7\,9\,10\}),$ $(\pi,\{4\,5\,6,3\,2\,4\,5\}),$ $(\pi,\{3\,2\,4\,5,8\,7\,9\,10\}),$ $(\pi,\{4\,5\,6\}),$ $(\pi,\{3\,2\,4\,5\}),$ $(\pi,\{8\,7\,9\,10\})$ and $(\pi,\emptyset)$ are all the possible $T$-marked words on $\pi.$  

A $T$-marked word $(\pi,M)$ is said to be the \emph{concatenation} of two (non-empty) $T$-marked words $(\alpha,M_1)$ and $(\beta,M_2)$ if $\pi=\alpha \beta$ and $M=M_1\cup M_2.$

A $T$-\emph{cluster} is a $T$-marked word $(\pi,M)$ that is not
the concatenation of two non-empty marked words and whose underlying word is a permutation of length at least two.

As an example, consider $T=\{2134,123\}$ and $\pi=2\,1\,3\,4\,5\,6\,7.$
The two $T$-marked words represented below

\begin{displaymath}\begin{tikzpicture}[   txt/.style = {text height=2ex, text depth=0.25ex}]
  \node [txt]  (0) at (0,0) {2};
  \node [txt]  (0) at (1/2,0) {1};
  \node [txt]  (0) at (2/2,0) {3};
  \node [txt]  (0) at (3/2,0) {4};
  \node [txt]  (0) at (4/2,0) {5};
  \node [txt]  (0) at (5/2,0) {6};
  \node [txt]  (0) at (6/2,0) {7};
\draw [decorate,decoration={brace,mirror,amplitude=5pt}]
    (-0.2,-0.2) -- (1.7,-0.2) ;
    \draw [decorate,decoration={brace,amplitude=5pt}]
    (0.8,0.4) -- (2.2,0.4) ;
   \draw [decorate,decoration={brace,mirror,amplitude=5pt}]
    (1.3,-0.4) -- (2.7,-0.4) ;   
   \draw [decorate,decoration={brace,amplitude=5pt}]
    (1.8,0.2) -- (3.2,0.2) ;  
 \end{tikzpicture}\end{displaymath}

\begin{displaymath}\begin{tikzpicture}[   txt/.style = {text height=2ex, text depth=0.25ex}]
  \node [txt]  (0) at (0,0) {2};
  \node [txt]  (0) at (1/2,0) {1};
  \node [txt]  (0) at (2/2,0) {3};
  \node [txt]  (0) at (3/2,0) {4};
  \node [txt]  (0) at (4/2,0) {5};
  \node [txt]  (0) at (5/2,0) {6};
  \node [txt]  (0) at (6/2,0) {7};
  \node [txt]  (0) at (3.3,0) {,};
\draw [decorate,decoration={brace,mirror,amplitude=5pt}]
    (-0.2,-0.2) -- (1.7,-0.2) ;
    \draw [decorate,decoration={brace,amplitude=5pt}]
    (0.8,0.4) -- (2.2,0.4) ;
   \draw [decorate,decoration={brace,amplitude=5pt}]
    (1.8,0.2) -- (3.2,0.2) ;  
 \end{tikzpicture}\end{displaymath}

namely, $(\pi,\{2134,345,456,567\})$ and $(\pi,\{2134,345,567\}),$ are $T$-clusters, while

\begin{displaymath}\begin{tikzpicture}[   txt/.style = {text height=2ex, text depth=0.25ex}]
  \node [txt]  (0) at (0,0) {2};
  \node [txt]  (0) at (1/2,0) {1};
  \node [txt]  (0) at (2/2,0) {3};
  \node [txt]  (0) at (3/2,0) {4};
  \node [txt]  (0) at (4/2,0) {5};
  \node [txt]  (0) at (5/2,0) {6};
  \node [txt]  (0) at (6/2,0) {7};
\draw [decorate,decoration={brace,mirror,amplitude=5pt}]
    (-0.2,-0.2) -- (1.7,-0.2) ;
     \draw [decorate,decoration={brace,amplitude=5pt}]
    (1.8,0.2) -- (3.2,0.2) ;  
 \end{tikzpicture}\end{displaymath}

and

\begin{displaymath}\begin{tikzpicture}[   txt/.style = {text height=2ex, text depth=0.25ex}]
  \node [txt]  (0) at (0,0) {2};
  \node [txt]  (0) at (1/2,0) {1};
  \node [txt]  (0) at (2/2,0) {3};
  \node [txt]  (0) at (3/2,0) {4};
  \node [txt]  (0) at (4/2,0) {5};
  \node [txt]  (0) at (5/2,0) {6};
  \node [txt]  (0) at (6/2,0) {7};
    \node [txt]  (0) at (3.3,0) {,};
    \draw [decorate,decoration={brace,amplitude=5pt}]
    (0.8,0.4) -- (2.2,0.4) ;
   \draw [decorate,decoration={brace,mirror,amplitude=5pt}]
    (1.3,-0.4) -- (2.7,-0.4) ;   
   \draw [decorate,decoration={brace,amplitude=5pt}]
    (1.8,0.2) -- (3.2,0.2) ;  
 \end{tikzpicture}\end{displaymath}

namely, $(\pi,\{2134,567\})$ and $(\pi,\{345,456,567\}),$ are not $T$-clusters. On the other hand, $(\pi,\{2134,567\})$ is the concatenation of $(2134,\{2134\})$ and $(567,\{567\}).$

An essential tool for studying Hertzsprung patterns is adapting the notion of inflation of permutations (see \cite{ALBERT20051}) to the context of marked permutations, as in \cite{Claesson_Hertz}.
Given a permutation $\rho$ of length $m$ and non-empty marked permutations $(\alpha_1,M_1),\ldots,(\alpha_m,M_m),$ the \emph{H-inflation} of $\rho$ by $(\alpha_1,M_1),\ldots,(\alpha_m,M_m),$ denoted $\rho[(\alpha_1,M_1),\ldots,(\alpha_m,M_m)]$, is the marked permutation $(\sigma,M),$ such that $\sigma=\alpha_1'\ldots \alpha_m'$ 
and $M=M'_1\cup\ldots \cup M'_m$
where  $\alpha'_{\sigma^{-1}(i)}$ is obtained from $\alpha_{\sigma^{-1}(i)}$ 
by adding the constant 
$|\alpha_{\sigma^{-1}(1)}|+\ldots+|\alpha_{\sigma^{-1}(i-1)}|$ and $M'_i$ is obtained from $M_i$ by adding the same constant to all its elements.

As an example, $$312[(1,\emptyset),(21,\{2\}),(123,\{12,23\})]=   (621345,\{2,34,45\}),$$ as depicted in the following figure, where the marks are represented by blue lines. 

\usetikzlibrary{arrows.meta,positioning}
\begin{figure}[h!]
\centering
\begin{tikzpicture}[scale=0.8]
    \draw[->] (0,0) -- (7,0) node[right] {$i$};
    \draw[->] (0,0) -- (0,7) node[above] {$\pi(i)$};

    \draw[red!30, fill=red!10, thick] (0.5,5.5) rectangle (1.5,6.5);
    \draw[blue!30, fill=blue!10, thick] (1.5,0.5) rectangle (3.5,2.5);
    \draw[green!30, fill=green!10, thick] (3.5,2.5) rectangle (6.5,5.5);

    \foreach \x/\y in {1/6, 2/2, 3/1, 4/3, 5/4, 6/5} {
        \filldraw[black] (\x, \y) circle (2pt);
    }

    \draw[blue, thick] (2,2) circle [radius=0.25];

    \draw[blue, thick, rotate around={45:(4.5,3.5)}] (4.5,3.5) ellipse (0.9 and 0.4);

    \draw[blue, thick, rotate around={45:(5.5,4.5)}] (5.5,4.5) ellipse (0.9 and 0.4);

\end{tikzpicture}
 \caption{Graphical representation of the inflation \(312[(1,\emptyset),(21,\{2\}),(123,\{12,23\})] = (621345,\{2,34,45\})\), with blocks highlighted and relations marked.}
\end{figure}

The following theorem (Theorem 2.1 in \cite{Claesson_Hertz}) ensures that every marked permutation can be written in a unique way as an inflation of clusters.

\begin{thm}\label{unique}
Let $T$ be an antichain of H-patterns. Any $T$-marked permutation can be uniquely written as $\rho[(\alpha_1,M_1),\ldots,(\alpha_m,M_m)],$ where $\rho \in S_m$ and each $(\alpha_i,M_i)$ is either the pair $(1,\emptyset),$ or a $T$-cluster.
\end{thm}

\section{General results}\label{general_results}

First of all we adapt the notions introduced in the previous section to the case of involutions. 
Let $\mathcal I_n$ be the set of involutions of length $n$ and $\mathcal I:=\cup_{n\geq 0}\mathcal I_n.$

A set $T$ of H-patterns is said to be \emph{self-inverse} if $\tau \in T$ whenever $\tau^{-1}\in T.$
Notice that, if the involution $\pi=\pi_1\pi_1\ldots \pi_n$ contains an occurrence $\pi_j\pi_{j+1}\ldots \pi_{j+k-1}$ of the H-pattern $\tau,$ then the symbols $j,j+1,\ldots, j+k-1$ form an occurrence of the H-pattern $\tau^{-1}$ in $\pi.$ Such occurrence will be called the \emph{sibling} of the occurrence $\pi_j\pi_{j+1}\ldots \pi_{j+k-1}.$

For example, the involution $\pi=10\, 8\,9\,7\,5\,6\,4\,2\,3\,1$ contains an occurrence of the H-pattern 231, namely, 897, whose sibling occurrence is 423, which is an occurrence of $312=231^{-1}.$

Let $T$ be a self-inverse set of H-patterns. A $T$-marked permutation  $(\pi,M)$ will be said a $T$-\emph{marked involution} if $\pi$ is an involution and $M$ contains an occurrence whenever it contains its sibling. 

We are now in position to state the analogue of Theorem \ref{unique} for involutions. 
\begin{thm}\label{unique_inv}
Let $T$ be an antichain of self-inverse H-patterns. Any $T$-marked involution can be written in a unique way as $\rho[(\alpha_1,M_1),\ldots,(\alpha_m,M_m)],$ where
\begin{itemize}
    \item $\rho$ is an involution in $\mathcal I_m$,
    \item each $(\alpha_i,M_i)$ is either the pair $(1,\emptyset),$ or a $T$-cluster, $1\leq i\leq m,$ 
    \item $\alpha_i=\alpha_{\rho_i}^{-1},$ $1\leq i\leq m,$ and
    \item $M_{\rho_i}$ contains precisely the siblings of the occurrences in $M_i,$ $1\leq i\leq m.$
\end{itemize}
\end{thm}

\proof 
Let $\pi$ be any $T$-marked involution. By Theorem \ref{unique}
there is a unique way to write $\pi$ as 

\noindent 
$\rho[(\alpha_1,M_1),\ldots,(\alpha_m,M_m)],$ where each $(\alpha_i,M_i)$ is either the pair $(1,\emptyset),$ or a $T$-cluster. Since $\pi$ is an involution, $\rho$ must be an involution itself. Moreover, when we inflate a fixed point $\rho_i$ of $\rho$ with the permutation $\alpha_i,$ in order to get an involution, $\alpha_i$ must be an involution. Similarly, if $(\rho_i,\rho_j)$ is a cycle of $\rho,$ if we inflate $\rho_i$ by $\alpha_i,$ $\rho_j$ must be inflated by $\alpha_i^{-1}.$ In this last case, by the definition of marked involution, the marked occurrences of $\alpha_i$ must be precisely the siblings of the marked occurrences of $\alpha^{-1}_i.$
\endproof

For example, let $T=\{231,312\}$ and let $$(\pi,M)=(10\, 8\,9\,7\,5\,6\,4\,2\,3\,1,\{10\,8\,9,7\,5\,6,5\,6\,4,2\,3\,1\})$$ be a $T$-marked involution. Then $(\pi,M)=\rho[(\alpha_1,M_1),(\alpha_2,M_2),(\alpha_3,M_3)],$ where
\begin{itemize}
\item $\rho=321,$
\item $(\alpha_1,M_1)=(312,\{312\}),$
\item $(\alpha_2,M_2)=(4231,\{423,231\}),$ and
\item $(\alpha_3,M_3)=(231,\{231\}),$
\end{itemize}
see figure below. 

\begin{figure}[h!]
\centering
\begin{tikzpicture}[scale=0.7]
    \draw[->] (0,0) -- (11,0) node[right] {$i$};
    \draw[->] (0,0) -- (0,11) node[above] {$\pi(i)$};

    \draw[red!30, fill=red!10, thick] (0.5,7.5) rectangle (3.5,10.5);   
    \draw[green!30, fill=green!10, thick] (3.5,3.5) rectangle (7.5,7.5); 
    \draw[blue!30, fill=blue!10, thick] (7.5,0.5) rectangle (10.5,3.5);  

    \foreach \x/\y in {
        1/10, 2/8, 3/9,
        4/7, 5/5, 6/6, 7/4,
        8/2, 9/3, 10/1
    } {
        \filldraw[black] (\x,\y) circle (2pt);
    }

    \draw[blue, thick, rounded corners=5pt]
        (0.7,7.7) rectangle (3.3,10.3);

    \draw[blue, thick, rounded corners=5pt]
        (3.7,4.7) rectangle (6.3,7.3);

    \draw[blue, thick, rounded corners=5pt]
        (4.7,3.7) rectangle (7.3,6.3);

    \draw[blue, thick, rounded corners=5pt]
        (7.7,0.7) rectangle (10.3,3.3);

\end{tikzpicture}
\caption{
A $T$-marked involution with $T=\{231,312\}.$
}
\end{figure}

In the following, we will denote by $\mathcal{C}_T$ the set of $T$-clusters and by $\mathcal{CI}_T$ the set of involutory $T$-clusters, namely, $T$-clusters whose underlying permutation is an involution. 

Let $T$ be a self-inverse set of H-patterns. Decompose $T$ as $T=T_I\,\dot{\cup}\, T_L\dot{\cup}\, T_L^{-1},$ where $T_I=T\cap \mathcal{I},$ and the set $\{T_L,T_L^{-1}\}$ is a  partition of the set of  $T\setminus T_I$ into two sets such that the second contains precisely the inverses of the first. 
Let $\pi$ be an involution and $\tau$ an H-pattern in $T.$ 
The number of occurrences of $\tau$ in $\pi$ coinciding with its sibling  will be denoted by $\tau^{sib}(\pi),$ while the number of occurrences of $\tau$ in $\pi$ not coinciding with its sibling  will be denoted by ${\tau^{nsib}(\pi)}.$ This last number must be even by definition. 

Now we introduce some generating functions useful in the following.
Let $\fp(\pi)$ and $\cyc(\pi)$ denote the number of \emph{fixed points} and number of \emph{cycles of length two} of $\pi,$ respectively. 

Let
$$I(t,y)=\sum_{\pi \in \mathcal I}t^{fp(\pi)}y^{cyc(\pi)}$$ be the generating function of involutions counted by number of fixed points and cycles.
It is well known that $I(t,y)$ can be expressed as a formal continued fraction as $$I(t,y)= \cfrac{1}{1-t-\cfrac{y}{1-t-\cfrac{2y}{1-t-\cfrac{3y}{\ldots}}}}.$$

This formula follows directly from the results in \cite{flaj}, where the author derives a more general expression for a generating function involving labeled Motzkin paths. The above generating function can be recovered as a suitable specialization, by exploiting a bijection between involutions and a specific subfamily of labeled Motzkin paths (see \textit{e.g.} \cite{Ba6}).

In order to define the generating function for $T$-clusters,
we need lists of variables $\{a_i\}_{i=1,\ldots ,r},$

\noindent
$\{b_j\}_{j=1,\ldots,s},$ and $\{c_k\}_{k=1,\ldots,s},$ that take into account the occurrences of H-patterns in $T_I=\{\tau_1,\ldots,\tau_r\},$  $T_L=\{\sigma_1,\ldots,\sigma_s\},$ and $T_L^{-1},$ respectively. 
Let 
\begin{equation*}C_T(x,a_1,\ldots a_{r},b_1,\ldots , b_{s},c_1,\ldots , c_{s})= \sum_{\pi\in \mathcal{C}_T}x^{|\pi|}\prod_{i=1}^r a_i^{{\tau_i(\pi)}}\prod_{j=1}^s b_j^{{\sigma_j(\pi)}} \prod_{k=1}^s c_k^{{\sigma^{-1}_k(\pi)}}\end{equation*}
be the \emph{cluster generating function}. 

We define also the generating function of involutory $T$-clusters. To this aim, we consider lists of variables $\{u_i\}_{i=1,\ldots,r},$ $\{\overline{u}_j\}_{j=1,\ldots,r},$ and $\{w_k\}_{k=1,\ldots,s},$  that take into account the occurrences of H-patterns in $T_I$ coinciding with their sibling, the occurrences of H-patterns in $T_I$ not coinciding with their sibling, and the occurrences of H-patterns in $T_L,$ respectively. 

Notice also that for every occurrence of $\sigma_i \in T_L$ in an involution there is precisely one occurrence of $\sigma_i^{-1}\in T_L^{-1}$ (its sibling occurrence), hence, when we deal with patterns contained in involutions, we can count only the occurrences of the former.  For this reason we define the the \emph{involutory cluster generating function} as follows.

\begin{equation*} CI_T(x,t,u_1,\ldots u_{r},\overline{u}_1,\ldots , \overline{u}_{r},w_1,\ldots , w_{s})= \sum_{\pi\in \mathcal{CI}_T}x^{|\pi|}t^{\fp(\pi)}\prod_{i=1}^r u_i^{{\tau^{sib}_i(\pi)}}\prod_{j=1}^r \overline{u}_j^{{\tau^{nsib}_j(\pi)}} \prod_{k=1}^s w_k^{{\sigma_k(\pi)}}.\end{equation*}

Finally, let

\begin{equation*} F_T(x,t,u_1,\ldots u_{r},\overline{u}_1,\ldots , \overline{u}_{r},w_1,\ldots , w_{s})= \sum_{\pi\in \mathcal I}x^{|\pi|}t^{\fp(\pi)}\prod_{i=1}^r u_i^{{\tau^{sib}_i(\pi)}}\prod_{j=1}^r \overline{u}_j^{{\tau^{nsib}_j(\pi)}} \prod_{k=1}^s w_k^{{\sigma_k(\pi)}}\end{equation*}
be the \emph{generating function for involutions} counted by length, number of fixed points and occurrences of H-patterns in $T.$

We are now in position to state our main result.

\begin{thm}\label{main_th}
Let $T$ be an antichain of self-inverse H-patterns. Then we have

\begin{multline*}F_T(x,t,u_1,\ldots u_{r},\overline{u}_1,\ldots , \overline{u}_{r},w_1,\ldots , w_{s})=\\ I\biggl(xt+CI_T(x,t,u_1-1,\ldots u_{r}-1,\sqrt{\overline{u}_1^2-1},\ldots , \sqrt{\overline{u}_{r}^2-1},w_1-1,\ldots , w_{s}-1),\\ x^2+C_T(x^2,{\overline{u}_1^2-1},\ldots , {\overline{u}_{r}^2-1},w_1-1,\ldots , w_{s}-1,w_1-1,\ldots , w_{s}-1)\biggr). \end{multline*}

\end{thm}
\proof
First of all, notice that, by definition of marked involution, the generating function  

$$MI_T(x,t,u_1,\ldots u_{r},\overline{u}_1,\ldots , \overline{u}_{r},w_1,\ldots , w_{s})$$ that
counts marked involutions by length ($x$), number of fixed points ($t$), number of marked occurrences of $\tau_i$ coinciding with their  sibling ($u_i$),  number of marked occurrences of $\tau_j$ not coinciding with their  sibling ($\overline{u}_j$), and  number of marked occurrences of $\sigma_k$ ($w_k$), can be written as

\begin{multline*}
MI_T(x,t,u_1,\ldots u_{r},\overline{u}_1,\ldots , \overline{u}_{r},w_1,\ldots , w_{s})=\\
\sum_{\pi\in \mathcal{I}}x^{|\pi|}t^{\fp(\pi)}\prod_{i=1}^r (1+u_i)^{{\tau^{sib}_i(\pi)}}\prod_{j=1}^r (1+\overline{u}_j^2)^{{\tau^{nsib}_j(\pi)}/2} \prod_{k=1}^s (1+w_k)^{{\sigma_k(\pi)}}.
\end{multline*}
In fact, the first product depends on the fact that, for every occurrence of $\tau_i$ coinciding with its sibling, we can choose either to mark it or not (the same holds for the third product). As for the second product, notice that every occurrence of $\tau_i$ not coinciding with its sibling yields a further occurrence of the same type. These two occurrences can be either both marked or both unmarked. 
Now, by Theorem \ref{unique_inv}, any marked involution $\pi$ can be uniquely written as 
 $\rho[(\alpha_1,M_1),\ldots,(\alpha_m,M_m)],$ where
\begin{itemize}
    \item $\rho$ is an involution,
    \item each $(\alpha_i,M_i)$ is either the pair $(1,\emptyset),$ or a $T$-cluster, 
    \item $\alpha_i=\alpha_{\rho_i}^{-1},$ and
    \item $M_{\rho_i}$ contains precisely the siblings of the occurrences of $M_i.$
\end{itemize}

For every $i=1,\ldots,m,$ consider the inflation of the entry $\rho_i$ by $(\alpha_i,M_i).$ Suppose that $(\alpha_i,M_i)$ is a $T$-cluster. 

We analyze two cases.

\begin{itemize}
    \item If $\rho_i$ is a fixed point and hence $\alpha_i$ is an involution, then
    \begin{itemize}
        \item any marked occurrence of a pattern $\tau\in T_I$ in $\alpha_i$ coinciding with its sibling turns into a marked occurrence of $\tau$ in $\pi$ coinciding with its sibling, 
        \item any marked occurrence of $\tau\in T_I$ in $\alpha_i$ not coinciding with its sibling turns into a marked occurrence of $\tau$ in $\pi$ not coinciding with its sibling.   Notice that in this case there is always the corresponding sibling occurrence both in $\alpha_i$ and in $\pi,$
        \item any marked occurrence of a pattern $\sigma\in T_L$ in $\alpha_i$ not coinciding with its sibling  turns into a marked occurrence of $\sigma$ in $\pi$ not coinciding with its sibling. Notice that also in this case there is always the corresponding sibling occurrence both in $\alpha_i$ and in $\pi.$
    \end{itemize}

\item If $\rho_i$ is not a fixed point of $\pi,$ we must inflate the two elements of the cycle $(\rho_i,\rho_j)$  by  $\alpha_i$ and $\alpha^{-1}_i,$ respectively.
In this case,
    \begin{itemize}
    \item any marked occurrence of a pattern $\tau\in T_I$ in $\alpha_i$ turns into a pair of distinct sibling marked occurrences of $\tau$ in $\pi,$ and
    \item any marked occurrence of a pattern $\sigma\in T_L$ in $\alpha_i$ turns into a marked occurrence of $\sigma$ in $\pi.$  Notice that in this case there is always the corresponding sibling occurrence both in $\alpha_i$ and in $\pi.$
    \end{itemize}
    
\end{itemize}

As a consequence,  

\begin{multline*}MI_T(x,t,u_1,\ldots u_{r},\overline{u}_1,\ldots , \overline{u}_{r},w_1,\ldots , w_{s})=\\ I\biggl( xt+CI_T(x,t,u_1,\ldots u_{r},\overline{u}_1,\ldots , \overline{u}_{r},w_1,\ldots , w_{s}),\\ x^2+C_T(x^2,\overline{u}_1^2,\ldots , \overline{u}_{r}^2,w_1,\ldots , w_{s},w_1,\ldots , w_{s})\biggr). \end{multline*}

This concludes the proof. 
\endproof

\section{Patterns 12 and 21}

As a first example of application of Theorem \ref{main_th}, we consider the set of H-patterns $T=\{12,21\}.$ 
Notice that in this case $T=T_I,$ hence $$\mathcal{C}=\mathcal{CI}$$
and such a set consists of all the clusters of the form $$(1\,2\,\ldots m,\,\{1\,2,2\,3,\ldots, m\hspace{-3pt}-\hspace{-3pt}1\,m\})$$ with $m\geq 2,$ and $$(n\,n\hspace{-3pt}-\hspace{-3pt}1\ldots 1,\,\{n\,n\hspace{-3pt}-\hspace{-3pt}1,\ldots ,3\,2, 2\,1\})$$
with $n\geq 2.$ Notice that the first one of these clusters contains $m$ fixed points, while the second one contains one fixed point if $n$ is odd and no fixed points otherwise. 

From these considerations it follows immediately that \\
$C(x,a_1,a_2)=CI(x,1,a_1,a_2,a_1,a_2)$ and that 
\begin{align*}
CI(x,t,u_1,u_2,\overline{u}_1,\overline{u}_2) &=\sum_{m\geq 2}x^m t^m u_1^{m-1}+\sum_{n\geq 1}(x^{2n}u_2\overline{u}_2^{2n-2}+x^{2n+1}t\overline{u}_2^{2n})\\
&=\frac{x^2t^2u_1}{1-xtu_1}+\frac{(u_2+xt\overline{u}_2^2)x^2}{1-x^2\overline{u}_2^2}.
\end{align*}

As a consequence of Theorem \ref{main_th} we have 

$$F(x,t,u_1,u_2,\overline{u}_1,\overline{u}_2)= \cfrac{1}{1-A-\cfrac{B}{1-A-\cfrac{2B}{1-A-\cfrac{3B}{\ldots}}}}$$
where $$A=xt+\frac{x^2t^2(u_1-1)}{1-xt(u_1-1)}+\frac{(u_2-1+xt(\overline{u}_2^2-1))x^2}{1-x^2(\overline{u}_2^2-1)},$$ and

\begin{align*}B &=x^2+\frac{x^4(\overline{u}_1^2-1)}{1-x^2(\overline{u}_1^2-1)}+\frac{(\overline{u}_2^2-1+x^2(\overline{u}_2^2-1)^2)x^4}{1-x^4(\overline{u}_2^2-1)^2}\\
&=x^2+\frac{x^4(\overline{u}_1^2-1)}{1-x^2(\overline{u}_1^2-1)}+\frac{x^4(\overline{u}_2^2-1)}{1-x^2(\overline{u}_2^2-1)}.
\end{align*}
The first terms of the generating functions $F$ are 
$$F(x,t,u_1,u_2,\overline{u}_1,\overline{u}_2)=1+tx+(t^2u_1 + u_2)x^2+ (t^3u_1^2 + t \overline{u}_2^2 + 2tu_2)x^3+$$
$$+( t^4u_1^3 + 2t^2u_1u_2 + 2t^2\overline{u}_2^2 + t^2u_1 + t^2u_2 + u_2\overline{u}_2^2 + u_2^2 + \overline{u}_1^2)x^4+\ldots$$

This result provides the solution of some enumerative problems by suitable specializations.

More precisely, setting $t=1$ and $u_1=u_2=\overline{u}_1=\overline{u}_2=0,$ we solve the Hertzprung problem (see \cite{Hertzprung}) over the set of involutions, \textit{i.e.}, the problem of determining the number of involutions such that two consecutive entries never differ by one. The first terms of the generating function in this case are $$F(x,1,0,0,0,0)=1+x+2x^5+8x^6+22x^{7}+74x^8+256x^9+\ldots$$

Moreover, setting $t=1,$ $u_1=\overline{u}_1=0$ and $u_2=\overline{u}_2=1,$ we enumerate \emph{irreducible involutions}, \textit{i.e.}, involutions $\pi$ such that $\pi(i+1)-\pi(i)\neq 1$ for every $i$ (see \cite{Baril} and seq. \seqnum{A278024} in \cite{Sl}).
Analogously, setting $t=0,$ $u_1=\overline{u}_1=0$ and $u_2=\overline{u}_2=1,$ we enumerate \emph{irreducible involutions} without fixed points (seq. \seqnum{A165968} in \cite{Sl}).

Finally, setting $t=0,$ $u_1=\overline{u}_1=\overline{u}_2=1$ we get the sequence of perfect matchings (involutions without fixed points) counted by occurrences of short-pairs, \textit{i.e.}, arcs connecting two consecutive elements (seq. \seqnum{A079267} in \cite{Sl}).

\section{Patterns $123$ and $321$}

Now we consider the set of H-patterns $T=\{123,\, 321\}.$
Also in this case $T=T_I,$   $\mathcal{C}=\mathcal{CI}$
and such a set contains all the clusters whose underlying permutations are$$1\,2\,\ldots m$$ with $m\geq 3,$ and $$n\,n\hspace{-3pt}-\hspace{-3pt}1\ldots 1$$ with $n\geq 3.$ 

Now we deduce an expression for the involutory cluster generating function $CI(x,t,u_1,u_2,\overline{u}_1,\overline{u}_2),$ where $u_1$ counts the occurrences of $123$ coinciding with their siblings, $u_2$ counts the occurrences of $321$ coinciding with their siblings, $\overline{u}_1$ counts the occurrences of $123$ not coinciding with their siblings, and $\overline{u}_2$ counts the occurrences of $321$ not coinciding with their siblings.
The generating function $CI$ is the sum of three terms $K_1,$ $K_2,$ and $K_3$ corresponding to the following cases.   

\begin{itemize}
\item Consider the clusters whose underlying permutation is $1\,2\,\ldots m.$ A cluster of this form contains $m$ fixed points and every occurrence of $123$ in it coincides with its sibling. 
Let $$K_1(x,t,u_1)=\sum_{m,j}a_{m,j}x^mt^mu_1^j,$$ where $a_{m,j}$ is the number of clusters over the permutation $12\ldots m$ with $j$ marked occurrences of $123.$ 

Let $\pi$ be a cluster of this kind, with $m\geq 5.$ The last marked occurrence in $\pi$ is $m\hspace{-3pt}-\hspace{-3pt}2\; m\hspace{-3pt}-\hspace{-3pt}1\; m.$ The second last one is either  $m\hspace{-3pt}-\hspace{-3pt}3\; m\hspace{-3pt}-\hspace{-3pt}2\; m\hspace{-3pt}-\hspace{-3pt}1,$ or $m\hspace{-3pt}-\hspace{-3pt}4\; m\hspace{-3pt}-\hspace{-3pt}3\; m\hspace{-3pt}-\hspace{-3pt}2.$ In the first case, if we remove the last symbol and the last marking, we get a cluster over $12\ldots m\hspace{-3pt}-\hspace{-3pt}1.$ In the second case, if we remove the last two symbols and the last marking, we get a cluster over $12\ldots m\hspace{-3pt}-\hspace{-3pt}2.$ Hence $$a_{m,j}=a_{m-1,j-1}+a_{m-2,j-1}.$$

This recurrence implies immediately 
$$K_1(x,t,u_1)=\frac{x^3t^3u_1}{1-u_1xt-u_1x^2t^2}.$$

\item Now we consider the clusters whose underlying permutation is $n\,n\hspace{-3pt}-\hspace{-3pt}1\ldots 1$ and with a marked occurrence of $321$ coinciding with its sibling (the central one). Such clusters have necessarily odd length, hence they have one fixed point. The contribution of such clusters to $CI$ is $$\frac{x^3tu_2}{1-\overline{u}_2^2x^2-\overline{u}_2^2x^4}.$$

Let $$K_2(x,t,u_2,\overline{u}_2)=\sum_{n,j}b_{n,j}x^ntu_2\overline{u}_2^j,$$ where $b_{n,j}$ is the number of clusters over the permutation $n\,n\hspace{-3pt}-\hspace{-3pt}1\ldots 1$ ($n$ odd) with $j$ marked occurrences of $321$ not coinciding with their sibling and one occurrence (the central one) coinciding with its sibling.  
Let $\pi$ be a cluster of this kind, with $n\geq 7.$ The first marked occurrence in $\pi$ is $n\; n\hspace{-3pt}-\hspace{-3pt}1\; n\hspace{-3pt}-\hspace{-3pt}2,$ whose sibling is the last marked occurrence, \textit{i.e.}, $321.$ The second one is either  $n\hspace{-3pt}-\hspace{-3pt}1\; n\hspace{-3pt}-\hspace{-3pt}2\; n\hspace{-3pt}-\hspace{-3pt}3$ or $n\hspace{-3pt}-\hspace{-3pt}2\; n\hspace{-3pt}-\hspace{-3pt}3\; n\hspace{-3pt}-\hspace{-3pt}4.$ In the first case, if we remove the first and the last symbol and the first and the last marking, after normalization we get a cluster over $n\hspace{-3pt}-\hspace{-3pt}2\,\,n\hspace{-3pt}-\hspace{-3pt}3\,\ldots 1.$ In the second case, if we remove the first and the last two symbols and first and the last marking, after normalization we get a cluster over $n\hspace{-3pt}-\hspace{-3pt}4\,\,n\hspace{-3pt}-\hspace{-3pt}5\ldots 1.$ Hence $$b_{n,j}=b_{n-2,j-2}+b_{n-4,j-2}.$$

This recurrence implies immediately 
$$K_2(x,t,u_2,\overline{u}_2)=\frac{x^3tu_2}{1-\overline{u}_2^2x^2-\overline{u}_2^2x^4}.$$




 


\item Now we consider the clusters whose underlying permutation is $n\,n\hspace{-3pt}-\hspace{-3pt}1\ldots 1$ ($n$ even or odd) and without marked occurrences of $321$ coinciding with their sibling.
Proceeding as above we obtain that the contribution to $CI$ of such clusters is
$$K_3(x,t,\overline{u}_2)=\frac{x^4\overline{u}^2_2+x^5t\overline{u}^2_2}{1-\overline{u}_2^2x^2-\overline{u}_2^2x^4}.$$

\end{itemize}

Summing all the these contributions we get

$$CI=\frac{x^3t^3u_1}{1-xtu_1-x^2t^2u_1}+\frac{x^3tu_2+x^4\overline{u}_2^2+x^5t\overline{u}_2^2}{1-x^2\overline{u}_2^2-x^4\overline{u}_2^2}.$$

The generating function for arbitrary clusters will follow by the trivial observation that $C(x,a_1,a_2)=CI(x,1,a_1,a_2,a_1,a_2).$ 
Exploiting Theorem \ref{main_th} we get 
$$F(x,t,u_1,u_2,\overline{u}_1,\overline{u}_2)= \cfrac{1}{1-A-\cfrac{B}{1-A-\cfrac{2B}{1-A-\cfrac{3B}{\ldots}}}}$$
where $$A=xt+CI(x,t,u_1-1,u_2-1,\sqrt{\overline{u}_1^2-1},\sqrt{\overline{u}_{2}^2-1})$$ and $$B=x^2+C(x^2,{\overline{u}_1^2-1},{\overline{u}_{2}^2-1}).$$

In this case we get $$F(x,t,u_1,u_2,\overline{u}_1,\overline{u}_2)=1+tx+(t^2 + 1)x^2+(t^3u_1 + tu_2 + 2t)x^3+$$
$$+(t^4u_1^2+2t^2u_2+4t^2+\overline{u}_2^2 +2)x^4+\ldots$$

\section{Patterns 231 and 312}

Now we consider the set $T=\{231,312\}.$ In this case $T=T_L\,\dot{\bigcup}\, T_L^{-1}\,,$ where $T_L=\{231\}.$
Set $\sigma=231.$ It is easily seen that 
\begin{multline*}\mathcal{C}=\{(231,\{231\}),(312,\{312\}),(4231,\{423,231\}),(45312,\{453,312\}),\\(564231,\{564,423,231\}),(645312,\{645,453,312\}),\\(7564231,\{756,564,423,231\}),(78645312,\{786,645,453,312\}),\ldots\}\end{multline*}
where the general element of this set is obtained by a sequence of alternating and overlapping occurrences of  $\sigma$ and $\sigma^{-1}.$
When an occurrence of $\sigma$ is followed by an occurrence of $\sigma^{-1},$ they overlap on one element, otherwise they overlap on two elements. Notice that, for every $k\geq 0,$ the set $\mathcal{C}$ contains 
\begin{itemize}
\item two elements of length $3k+3$. Such two elements are not involutions. One of them has $k+1$ occurrences of $\sigma$ and $k$ occurrences of $\sigma^{-1}$ and the other one has $k$ occurrences of $\sigma$ and $k+1$ occurrences of $\sigma^{-1}.$
\item A single element $\pi$ of length $6k+4.$ This element is an involution, $\fp(\pi)=2$ and $\sigma(\pi)=\sigma^{-1}(\pi)=2k+1.$
\item A single element $\pi$ of length $6k+5.$    $\pi$ is an involution, $\fp(\pi)=1$ and $\sigma(\pi)=\sigma^{-1}(\pi)=2k+1.$ 
\item A single element $\pi$ of length $6k+7.$  $\pi$ is an involution,
 $\fp(\pi)=1$ and $\sigma(\pi)=\sigma^{-1}(\pi)=2k+2.$
 \item A single element $\pi$ of length $6k+8.$  $\pi$ is an involution,
$\fp(\pi)=2$ and $\sigma(\pi)=\sigma^{-1}(\pi)=2k+2,$
\end{itemize}
where we recall that $\sigma(\pi)$ denotes the number of occurrences of the H-pattern $\sigma$ in $\pi.$ Hence, 
\begin{multline*}\mathcal{CI}=\{(4231,\{423,231\}),(45312,\{453,312\}),\\(7564231,\{756,564,423,231\}),(78645312,\{786,645,453,312\}),\ldots\}\end{multline*}

From the previous observations we can deduce the expression for the generating function of clusters 

\begin{align*}C(x,b,c)&=\sum_{k\geq 0}\left(x^{3k+3}b^{k+1}c^k+x^{3k+3}b^{k}c^{k+1}+x^{3k+4}b^{k+1}c^{k+1}+x^{3k+5}b^{k+1}c^{k+1}\right)\\
&=x^3\left(b+c+xbc+x^2bc\right)\sum_{k\geq 0}x^{3k}b^kc^k=\frac{x^3\left(b+c+xbc+x^2bc\right)}{1-x^3bc},
\end{align*}
and involutory clusters
\begin{align*}CI(x,t,w) &=\sum_{k\geq 0}\left(x^{6k+4}w^{2k+1}t^{2}+x^{6k+5}w^{2k+1}t+x^{6k+7}w^{2k+2}t+x^{6k+8}w^{2k+2}t^{2}\right)\\
&=x^4wt\left(t+x+x^3w+x^4wt\right)\sum_{k\geq 0}x^{6k}w^{2k}=\frac{x^4wt\left(t+x+x^3w+x^4wt\right)}{1-x^6w^2}.
\end{align*}
As a consequence of Theorem \ref{main_th} we have 
$$F(x,t,w)= \cfrac{1}{1-A-\cfrac{B}{1-A-\cfrac{2B}{1-A-\cfrac{3B}{\ldots}}}}$$
where $$A=xt+\frac{x^4(w-1)t\left(t+x+(w-1)(x^3+x^4t)\right)}{1-x^6(w-1)^2},$$ and $$B=x^2+\frac{x^6\left(2w-2+x^2(w-1)^2+x^4(w-1)^2\right)}{1-x^6(w-1)^2}.$$

Here we get $$F(x,t,w)=1+tx+(t^2 + 1)x^2+(t^3 + 3t)x^3+(t^4 + t^2w + 5t^2 +3)x^4+\ldots$$

\section{Patterns 132 and 213}

Finally we consider the set $T=\{132,213\}.$ In this case $T=T_I.$ It is easily seen that 
\begin{multline*}\mathcal{C}=\mathcal{CI}=\{(132,\{132\}),(213,\{213\}),(1324,\{132,324\}),(21354,\{213,354\}),\\(132465,\{132,324,465\}),(213546,\{213,354,546\}),\\(1324657,\{132,324,465,657\}),(21354687,\{213,354,546,687\}),\ldots\}\end{multline*}
where the general element of the set is obtained by a sequence of alternating and overlapping occurrences of  $132$ and $213.$
When an occurrence of $132$ is followed by an occurrence of $213$ they overlap on two elements, otherwise they overlap on one element. 
Every cluster is involutory and every occurrence of a pattern coincides with its sibling. 
Notice that, for every $k\geq 1,$ the set 
$\mathcal{CI}$ contains
\begin{itemize}
    \item two elements $\pi_1$ and $\pi_2$ of length $3k,$ with $\fp(\pi_1)=\fp(\pi_2)=k,$ $132(\pi_1)=213(\pi_2)=k,$ and $213(\pi_1)=132(\pi_2)=k-1,$
    \item one element $\pi$ of length $3k+1$ with $\fp(\pi)=k+1,$ and $132(\pi)=213(\pi)=k,$
    \item one element $\pi$ of length $3k+2$ with $\fp(\pi)=k$ and $132(\pi)=213(\pi)=k.$
\end{itemize}
From the previous observations, we can deduce the expression for the generating function of involutory clusters 
\begin{align*}CI(x,t,u_1,u_2)&\\ &=\sum_{k\geq 1}\left(x^{3k}t^k(u_1^ku_2^{k-1}+u_1^{k-1}u_2^{k})+x^{3k+1}t^{k+1}(u_1u_2)^k+x^{3k+2}t^{k}(u_1u_2)^k\right)\\
&=\frac{(u_1+u_2)x^3t+x^4t^2u_1u_2+x^5tu_1u_2}{1-x^3tu_1u_2}
\end{align*}
and $C(x,a_1,a_2)=CI(x,1,a_1,a_2).$ 
As a consequence of Theorem \ref{main_th} we have 
$$F(x,t,u_1,u_2,\overline{u}_1,\overline{u}_2)= \cfrac{1}{1-A-\cfrac{B}{1-A-\cfrac{2B}{1-A-\cfrac{3B}{\ldots}}}}$$
where $$A=xt+CI(x,t,u_1-1,u_2-1,\sqrt{\overline{u}_1^2-1},\sqrt{\overline{u}_{2}^2-1})$$ and $$B=x^2+C(x^2,{\overline{u}_1^2-1},{\overline{u}_{2}^2-1}).$$

Here we get 
$$F(x,t,u_1,u_2,\overline{u}_1,\overline{u}_2)=1+tx+(t^2 + 1)x^2+(t^3 + tu_1 + tu_2 + t)x^3+$$ $$+(t^4 + t^2u_1 + 3t^2 + t^2u_1u_2 + t^2u_2 + 3)x^4+\ldots$$

\section{Final remarks}

Suitable specializations of the generating functions obtained in the previous Sections allow us to determine the distribution of the occurrences of every H-pattern of length 2 or 3. As a byproduct, we can determine all Wilf-equivalences among patterns of length 2 or 3, namely, we can partition the set of H-patterns of the same length into equivalence classes where for every H-pattern in a class the number of involutions of length $n$ avoiding the pattern is the same, for all $n.$ Quite surprisingly, the H-patterns 12 and 21 are not Wilf-equivalent, and the only Wilf-equivalences among H-patterns of length 3 are trivial, namely, the pattern 132 is equivalent to 213 and the pattern 231 is equivalent to 312. In fact, 132=213$^{rc}$, 231=312$^{rc}$, and the set of involutions is closed under the $rc$-operator (we recall that, given a permutation $\sigma=\sigma(1) \ \sigma(2) \ \cdots\ \sigma(n)$, the permutation $\sigma^{rc}$ is defined by $\sigma^{rc}(i)=n+1-\sigma(n+1-i)$).

\acknowledgements
We would like to express our  gratitude to the anonymous referees, whose inspiring comments and careful reading enhanced the quality of this paper.

\nocite{*}
\bibliographystyle{abbrvnat}
\bibliography{biblio_hertz}
\label{sec:biblio}

\end{document}